\documentclass[aps,pra,onecolumn,superscriptaddress,groupedaddress]{revtex4}

\usepackage{graphicx,amsmath,amssymb,mathtools,xcolor}
\numberwithin{equation}{section}

\definecolor{deepred}{HTML}{CD2626}

\DeclarePairedDelimiter\norm{\lVert}{\rVert}%

\begin{document}
\title{Connectivity of Soft Random Geometric Graphs Over Annuli}
\author{Alexander P. Giles}
\email[]{Alexander.Giles@bristol.ac.uk}
\address{School of Mathematics, University of Bristol, University Walk, Bristol BS8 1TW, United Kingdom}
\author{Orestis Georgiou}
\email[]{Orestis.Georgiou@toshiba-trel.com}
\address{Toshiba Telecommunications Research Laboratory, 32 Queen Square, Bristol, BS1 4ND, United Kingdom}
\author{Carl P. Dettmann}
\email[]{Carl.Dettmann@bristol.ac.uk}
\address{School of Mathematics, University of Bristol, University Walk, Bristol BS8 1TW, United Kingdom}

\date{\today}

\begin{abstract}
Nodes are randomly distributed within an annulus 
(and then a shell) to form a point pattern of communication terminals
which are linked stochastically according to the Rayleigh fading
of radio-frequency data signals. We then present analytic formulas for the connection probability
of these spatially embedded graphs, describing the connectivity behaviour as a dense-network limit is approached. This 
extends recent work modelling \textit{ad hoc} networks in non-convex domains.
\end{abstract}

\maketitle

\section{Introduction}\label{sec:one}
Soft random geometric graphs \cite{penrose2013} are 
network structures \cite{newmanbook} consisting
of a set of nodes placed according to a point process in some domain $\mathcal{{V}}\subseteq\mathbb{R}^{d}$
mutually coupled with a probability dependent on their Euclidean separation.
Examples of their current application include modelling the collective behavior
of multi-robot swarms \cite{brecht2013}, disease surveillance
\cite{eubank2004}, electrical smart grid engineering \cite{amin2013} and particularly our focus, communication theory \cite{tsebook}, 
where random geometric graphs have recently been used to model \textit{ad hoc}
wireless networks \cite{mao2011,orestis2014,hagenggi2009,li2009,tohbook,halldorsson2014} sharing information
over communication channels which have a complex, stochastic impulse response \cite{cef2012,orestis2014,orestis2015}.

In an earlier form, random `hard' or `unit disk' \cite{clark1991} geometric graphs are formed
by picking a finite number of points from a $d$-dimensional Euclidean metric space $\mathcal{V}$
(such as the unit square) which are then joined whenever they lie within some critical distance of each other. These networks
take their structure from the underlying planar topology of the (usually bounded)
set $\mathcal{V}$ in which they live, distinguishing them from the non-spatial random graphs of Erd\H{o}s and R\'{e}nyi \cite{erdos19602}, and were 
introduced by Gilbert \cite{gilbert1961} at (what was then) AT\&T Bell Telephone Laboratories.

This deterministic connection can, however, be generalised to probabilistic (or `soft') connection
\cite{penrose2013,cef2012,mao2011} in order to model signal fading. Commonly known as the `random connection model', 
we now have a connection function $H\left(\norm{x-y}\right)$ giving the probability that links will form
between nodes $x,y \in \mathcal{V}$ of a certain Euclidean displacement $\norm{x-y}$. This is a 
(much) more realistic model than that of the unit disks.
In a band-limited world of wireless communications continuously pressed
for the theoretical advances that can enable 5G cellular performance, this is an important new flexibility in the
model.

Connectivity is a central focus of much of the research \cite{cef2012,orestis2013, mao2011}.
For example, in \cite{cef2012}
(using a cluster expansion technique from statistical physics), at high node density $\rho$ 
the connection probability of a soft random geometric graph formed within a bounded domain
$\mathcal{{V}}$ is approximated as (the complement of) the
probability that exactly one isolated node appears in an otherwise connected graph.
This is justified by a conjecture of Penrose \cite{penrose2013} (which can be proved under more
restrictive conditions than considered here), asserting
that the number of isolated nodes follows a Poisson distribution whose mean
quickly decays as $\rho\rightarrow\infty$, thus highlighting the
impact of the domain's enclosing boundary \cite{cef2012, mao2011,kog105} where node isolation is most common.

Internal boundaries, such as obstacles, cause similar problems, and we focus our efforts on how this particular aspect of the domain
effects the graph behaviour. We therefore extend recent work on connectivity
within non-convex domains \cite{orestis2013,almiron2013,bocus2013,dettmann2014,coon2014}
(such as those containing internal walls \cite{orestis2013} or a complex, fractal
boundary \cite{dettmann2014}), deriving analytic formulas for the connection probability
$P_{fc}$ of soft random geometric graphs formed within the annulus and spherical shell geometries,
quantifying how simple convex obstacles (of radii $r$) affect connectivity. 
Specifically, we consider the situation where nodes connect with
a probability decaying exponentially with their mutual Euclidean separation 
(which is equivalent to `Rayleigh' fading, commonly found in models of signal propagation within cluttered,
urban environments).
	
This paper is structured as follows:
In section \ref{sec:two} we describe our model and define our graph-theoretic problem. In section \ref{sec:three} we evaluate the graph connection
probability for the annulus domain $\mathcal{{A}}$, 
incorporating both small and large circular obstacles as internal perimeters of radius $r$.
We then extend the annulus into its three-dimensional analogue known as the 
spherical shell $\mathcal{{S}}$ in section \ref{sec:four}. After discussing our results in section \ref{sec:five},
we conclude in section \ref{sec:six}.

\begin{figure}
	
\noindent \begin{centering}
\includegraphics[scale=0.20]{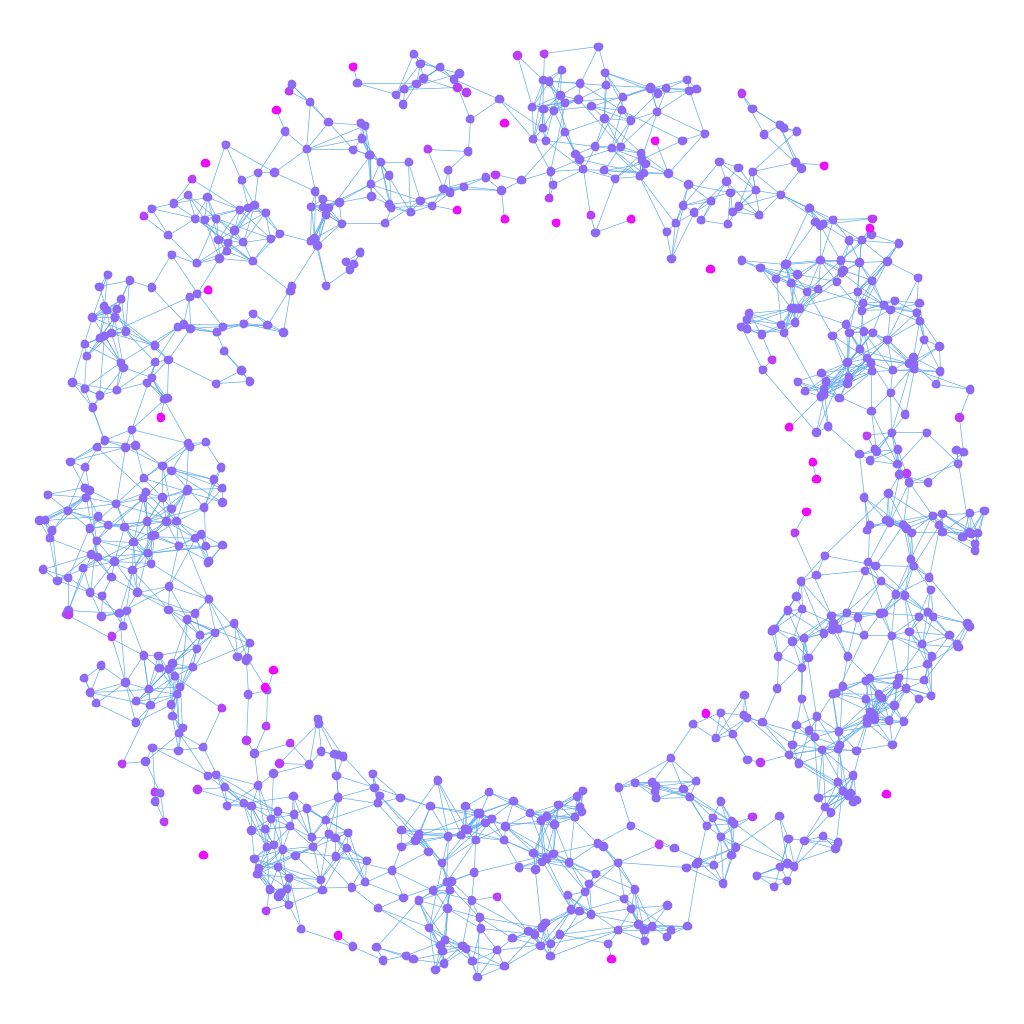}$\qquad\:$
\includegraphics[scale=0.20]{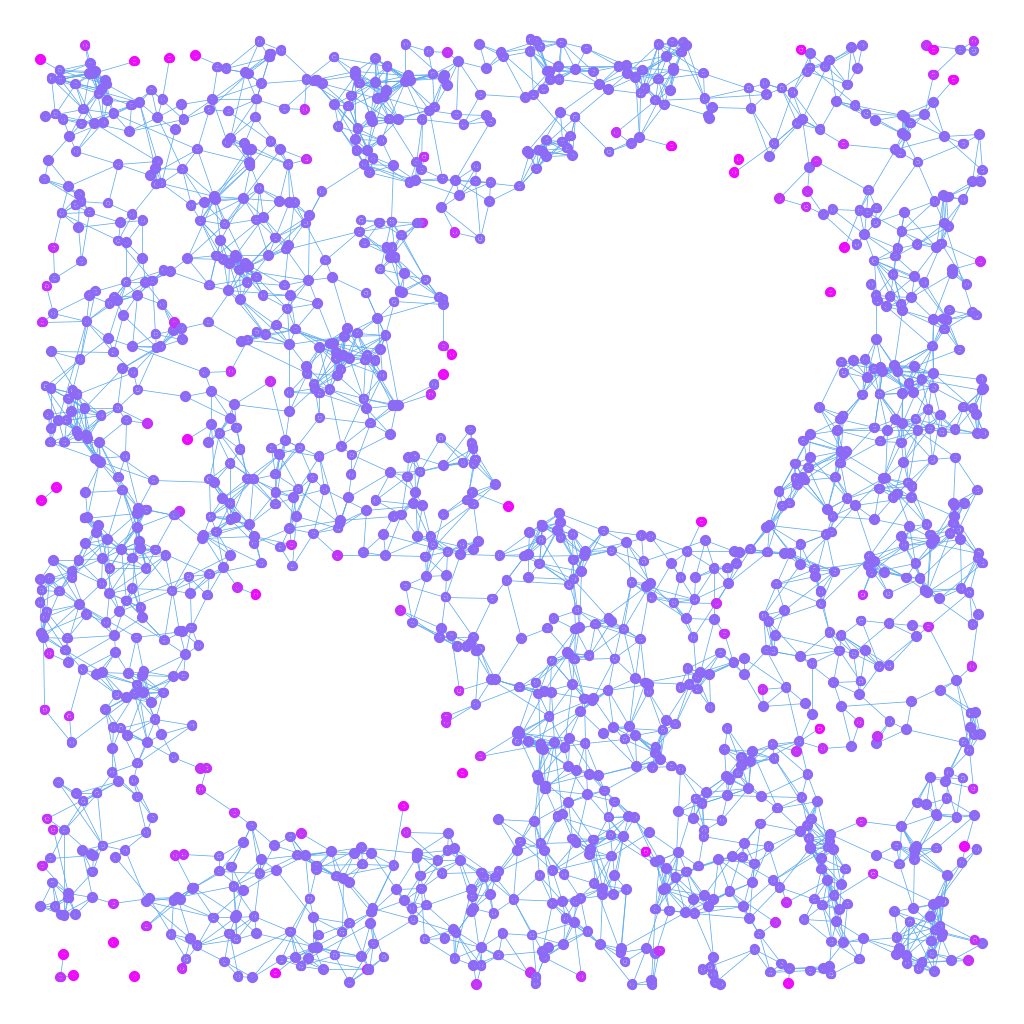}
\par\end{centering}

\caption{{\footnotesize (Colour online) Randomly constructed soft geometric graphs (picked from the
ensemble, with links formed based on the Rayleigh fading of electromagnetic data signals and nodes configured randomly) drawn inside $\mathcal{{A}}$ (left) and a square domain (right) containing two circular obstacles. Nodes with only a few connections (usually near the domain edges) are highlighted in purple, demonstrating the boundary effect phenomenon.}}
\label{fig:annulus}
\end{figure}

\section{Soft random geometric graphs bounded within non-convex geometries}\label{sec:two}

Let $\mathcal{V} \subseteq \mathbb{R}^{d}$ be a bounded region of volume $V$ associated with both the Lesbegue measure 
$\mathrm{d}x$ and the Euclidean metric $r_{xy} = \norm{x-y}$ for any $x,y \in \mathcal{V}$. We define
the characteristic function
\begin{equation}\label{e:neweq1}
	\chi \left(x,y\right) =
  \begin{cases}
	  1       & \quad \text{if } x + \lambda\left(y-x\right) \in \mathcal{V} \text{ for all } \lambda \in \left[0,1\right] \\
	  0  & \quad \text{otherwise} \\
  \end{cases}
\end{equation}
implying that $\mathcal{V}$ is a convex set whenever $\chi \left(x,y\right)=1$ for every pair $x,y \in \mathcal{V}$. 

Let $\mathcal{Y}$ be a Poisson point process of intensity $\rho $ on $\mathcal{V}$ with respect to Lesbegue measure $\textrm{d}x$ on $\mathbb{R}^d$. Construct a \textit{soft random geometric graph} in the following way: for every (unordered) pair of nodes $x, y \in \mathcal{Y}$, add an edge between $x$ and
$y$ (independently for each pair) with probability $\chi \left(x,y\right)H\left(r_{xy}\right)$, where the measurable function $H: \mathbb{R}^{d} \to \left[0,1\right]$ is called the \textit{connection function}. Should a node form $k$ links we say its has \textit{degree} $k$, and
call these linked nodes \textit{neighbours}.
\subsection{The degree distribution}
Given that $\mathcal{Y}$ contains a point at $x \in \mathcal{V}$, the remaining points of $\mathcal{Y}$ are
again distributed as a Poisson process of intensity $\rho$ on $\mathcal{V}$, by Palm theory
for the Poisson process: see e.g. $\S 1.7$ of Penrose's book \cite{penrosebook}. So consider a
vertex at some fixed $x \in \mathcal{V}$. We can determine its degree distribution by
looking at a marked Poisson point process $\mathcal{Y}^{\star}$ constructed by adorning
each $y \in \mathcal{Y}$ with an independent $U\left[0,1\right]$ random mark $u_y$, where $U\left[0,1\right]$
denotes the uniform distribution on the unit interval  $\left[0,1\right]$. Then by the
marking theorem for Poisson processes (see $\S 5.2$ of Kingman's book \cite{kingmanbook}),
\begin{eqnarray}
	\mathcal{Y}^{\star} = \left\{\left(y,u_y\right) : y \in \mathcal{Y}\right\}
\end{eqnarray}
is a Poisson point process on $\mathcal{V} \times \left[0,1\right]$ of intensity $\rho$ but now with respect
to Lebesgue measure on $\mathbb{R}^{d+1}$. The mark $u$ on each marked point $(y,u) \in \mathcal{Y}^{\star}$ is used to determine connectivity of $y$ to $x$: if $u \leq \chi(x,y)H(r_{xy})$, then
$y$ is joined to $x$ by an edge. The degree of $x$ is then distributed as
\begin{equation}\label{e:neweq2}
	k\left(x\right) = \sum_{\left(y,u\right) \in \mathcal{Y}^{\star}} \mathbf{1}\{u < \chi \left(x,y\right) H\left(r_{xy}\right)\}
\end{equation}
This representation for $k(x)$ as a sum over values of a measurable function $f$ of $(y,u) \in \mathcal{Y}^{\star}$ enables us to apply Campbell's theorem (see $\S 3.2$
of Kingman's book \cite{kingmanbook}) to obtain that $k(x)$ follows a Poisson distribution
with mean
\begin{equation}\label{e:neweq3}
	\mathbb{E}\left[k\left(x\right)\right] = \rho \int_{\mathcal{V} \times \left[0,1\right]} \mathbf{1}\{u < \chi \left(x,y\right) H\left(r_{xy}\right)\}\textrm{d}y\textrm{d}u =
	\rho \int_{\mathcal{V}}\chi \left(x,y\right) H\left(r_{xy}\right)\textrm{d}y
\end{equation}
In particular, the probability that $x$ has degree zero is
\begin{eqnarray}\label{e:neweq4}
	\exp\left(-\rho \int_{\mathcal{V}} \chi \left(x,y\right) H\left(r_{xy}\right) \mathrm{d}y\right)
\end{eqnarray}

Following Penrose's conjecture mentioned in Section \ref{sec:one} (and found as Theorem 2.1 in \cite{penrose2013}),
it is natural in light of Eq. \ref{e:neweq4} to
conjecture that as $\rho \to \infty$, the total number of isolated nodes
is well approximated by a Poisson distribution with mean
\begin{equation}\label{e:neweq5}
	\rho \int_{\mathcal{V}} \exp \left(-\rho \int_{\mathcal{V}} \chi \left(x,y\right) H\left(r_{xy}\right) \mathrm{d}y \right) \mathrm{d}x
\end{equation}

In this limit, the usual situation is that the obstacle to connectivity is the presence of
isolated nodes \cite{waltersreview}, and so another reasonable conjecture is that as $\rho \to \infty$

\begin{align} \label{e:neweq6}
	\mathbb{P}\left(\text{graph is connected}\right) &\approx \mathbb{P}\left(\text{no isolated nodes}\right) \nonumber  \\
													 &\approx \exp \left(-\rho \int_{\mathcal{V}} e^{-\rho \int_{\mathcal{V}} \chi \left(x,y\right) H\left(r_{xy}\right) \mathrm{d}y } \mathrm{d}x \right) \nonumber \\
				 &\approx 1 - \rho \int_{\mathcal{V}} e^{-\rho \int_{\mathcal{V}} \chi \left(x,y\right) H\left(r_{xy}\right) \mathrm{d}y } \mathrm{d}x
\end{align}

In the following article we evaluate this formula for various bounded, non-convex regions $\mathcal{V}$ in order to elucidate the specific effect of obstacles on high-density network connectivity. Our formulas
are `semi-rigorous' in that they are based on (at least) the above assumptions.
Rigorous proof of any formulas herein presented (in a
similar fashion to the work of Penrose) is deferred to a later study.

We also note that though we use the Poisson model for the point set $\mathcal{Y}$ throughout, our simulations (in Fig. \ref{fig:montecarlo}) consider only the \textit{binomial model}, where $N$ nodes are selected uniformly
at random from $\mathcal{V}$. The two models are closely related when $\rho=N/V$, in which case the Poisson process is locally
a good approximation to (what is called) the binomial point process of $N$ nodes \cite{kingmanbook}. 

\subsection{Dense networks}

It is important to highlight our scaling of density and volume. It is common to find asymptotic results for
connectivity in the literature (see e.g. \cite{gupta1998}), where points are drawn according to the usual
Poisson process of intensity $\rho \textrm{d}x$ inside a square $S_n$ of area $n$. Then, with $\rho$ fixed, 
one studies the limit $n \to \infty$. This would represent the thermodynamic limit if the critical connection range
$r$ were to remain independent of the geometry of the square, 
but it does not and instead we have $r\left(n\right)$ scaling in some way with
the geometry. 

Call graphs formed in this way $G\left(n,r\left(n\right)\right)$. One can then prove
that, given some supercritical connection range
\begin{eqnarray}
	\lim_{n \to \infty} \mathbb{P}\left(G\left(n,r\left(n\right)\right) \textit{is connected}\right) = 1
\end{eqnarray}
This critical scaling is known \cite{penrosebook}
to be of the order $r^2\left(n\right) \approx \log n$, beyond which we enter a phase where \textit{almost all} graphs connect rather than
\textit{almost all} graphs do not connect. One is interested in the point in the parameter space
at which this transition occurs.

But our random geometric graphs do not scale this way. The nodes connect
according to a continuous function of their Euclidean separation, and we asymptotically scale the geometry-independent
density of the point process $\rho \to \infty$. As has previously been noted \cite{cef2012,coon2012}, the
enclosing boundary of $\mathcal{V}$ becomes the dominant influence on connectivity as density increases. We study
a similar effect, but concerning convex obstacles, rather than the enclosing perimeter.

\subsection{Rayleigh fading means soft connection}

In `soft' graphs, links are formed between nodes with a certain probability (such as according
to a coin toss).  This allows us to model \textit{inter alia} data transmission in an environment 
where connection fails according
to an activity known as \textit{signal fading}, a simple example
of which is \textit{Rayleigh fading} \cite{sklar1997,tsebook}.

The information outage probability $P_{out}$ (quantifying how often
the signal's distortion causes the decoding error at the receiver to fall
below a threshold rate $\Upsilon$)
is given in the Rayleigh model by
\begin{eqnarray}\label{e:neweq8}
P_{out} & = & P\left[\log_{2}\left(1+\frac{P}{N_{0}}\hspace{0.5mm}\norm{h}^{2}\right)<\Upsilon\right]
\end{eqnarray} where $h$ is the channel gain and $P/N_0$ is the signal-to-noise
ratio. Since the channel gain is a complex Gaussian process, its magnitude is Rayleigh distributed and so
$\norm{h}^{2}$ is exponentially distributed
\begin{eqnarray}\label{e:neweq10}
H\left(r_{xy}\right) &=& 1 - P\left[ \norm{h}^{2} < \frac{ N_{0} \left(2^{\Upsilon}-1\right)}{P}\right] \nonumber \\
									& = &	1 - P\left[ \norm{h}^{2} < \beta r_{xy}^{\eta} \right] \nonumber \\
& = &   e^{-\beta r^{\eta}_{xy}}
\end{eqnarray}
since the signal to noise ratio is related to the propagation distance
$r_{xy}$ through $P = c r_{xy}^{\eta}$ where $\eta$ is the `path loss exponent'; we uniquely consider
$\eta = 2$ (free-space propagation). Thus $\beta=N_{0} \left(2^{\Upsilon}-1\right) / c$. 
We also use
\begin{eqnarray}
	r_{0} = \beta^{-1/\eta}
\end{eqnarray}
to signify the length scale over which nodes connect,
since the exponent $\beta r_{xy}^\eta>1$ whenever $r_{xy}>r_0$ (and
\textit{vice versa}).

\section{The annulus domain $\mathcal{{A}}$}\label{sec:three}

\begin{figure}
\noindent \begin{centering}
\includegraphics[scale=0.3]{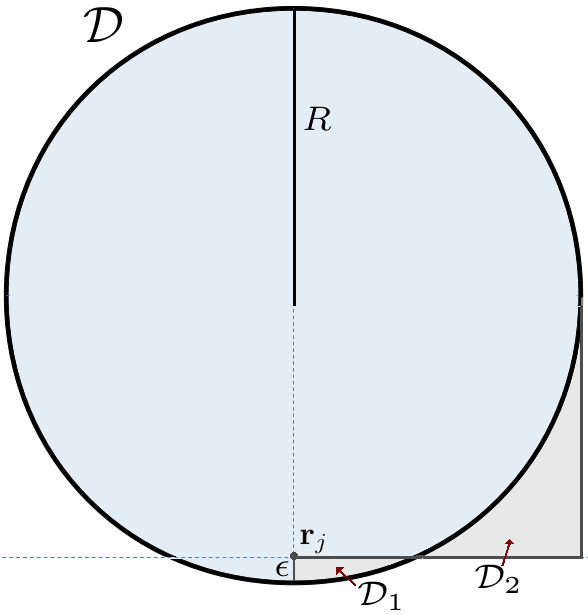}\qquad\includegraphics[scale=0.3]{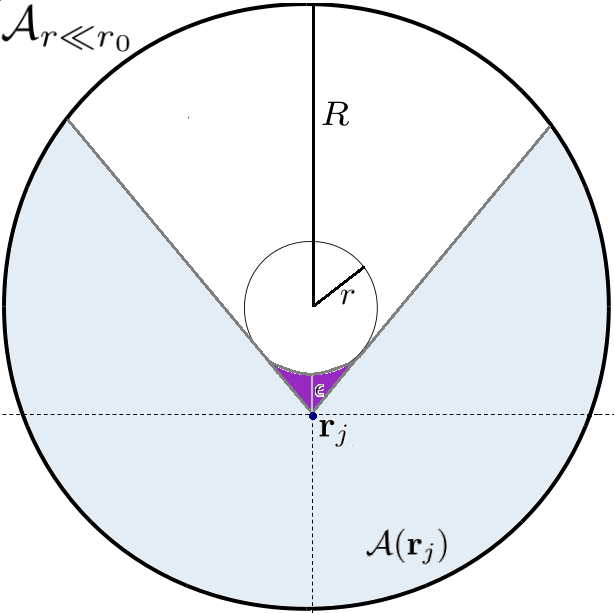}\qquad\includegraphics[scale=0.3]{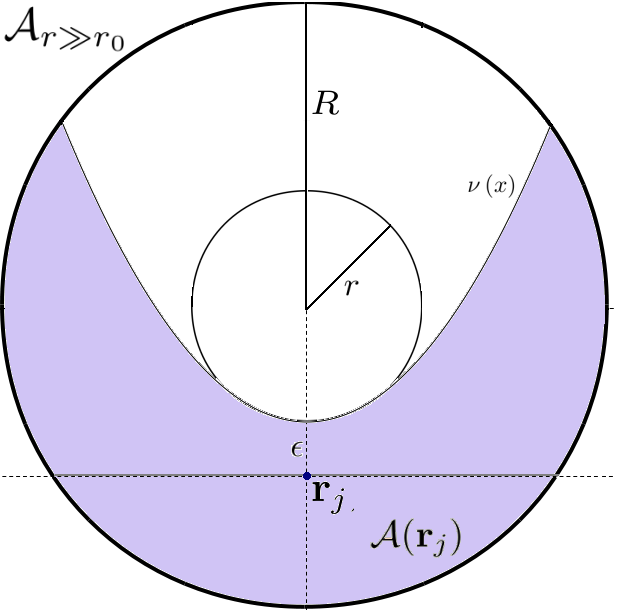}
\par\end{centering}
\caption{{\footnotesize (Colour online) A depiction of the integration regions used for the disk domain $\mathcal{{D}}$ (left panel) and annulus domain $\mathcal{{A}}$ with small obstruction (middle panel) and large obstruction (right panel), with the integration regions highlighted. The small, cone-like region in the middle domain
$\mathcal{A}$ is highlighted in purple.}\label{fig:trip}}
\end{figure}

In this section, we take $\mathcal{V}$ to be the 
annulus $\mathcal{{A}}$ of inner radius $r$ and outer radius $R \gg r_0$ (depicted in the left panel of Fig. \ref{fig:trip}). 
In order to simplify the necessary integrals, we define the `connectivity mass' at
$x \in \mathcal{A}$ 
\begin{eqnarray}
	\mathcal{M}\left(x\right) &=& \int_{\mathcal{A}} \chi\left(x,y\right) H\left(r_{xy}\right)\textrm{d}y \nonumber \\
															   &=&  \int_{\mathcal{A}\left(x\right)} H\left(r_{xy}\right)\textrm{d}y
\end{eqnarray}
(taken from the exponent in Eq. \ref{e:neweq6}). This is approximated within two obstacle-size regimes, the first
where $r\ll r_{0}$, and the second where $r\gg r_{0}$;
in each regime we can make some assumptions about the geometry of the region
$\mathcal{A}\left(x\right)$ visible to $x$, which yields tractable formulas for the connectivity mass in terms of
powers of the distance $\epsilon \in \left[0,R-r\right]$ from the obstacle perimeter. Given a slight correction
to a previous result in \cite{cef2012} on connectivity within a disk of radius $R$, we then have $P_{fc}$ in $\mathcal{A}$.

\subsection{No obstacles}

We first take the case where $r=0$ depicted in Fig \ref{fig:trip} (which is the disk $\mathcal{D}$). We quickly derive
an approximation to $P_{fc}$ in this limiting domain (which we later extend into the annulus $\mathcal{A}$).

We first have the connectivity mass a distance $\epsilon$ from the disk's
centre, given by
\begin{eqnarray} \label{eq:neweq16}
\mathcal{{M}}\left(\epsilon\right) & = & \frac{\pi}{2\beta}+2\left(\int_{\mathcal{{D}}_{1}}
e^{-\beta\left(x^{2}+y^{2}\right)}\textrm{d}y\textrm{d}x-\int_{\mathcal{{D}}_{2}}
e^{-\beta\left(x^{2}+y^{2}\right)}\textrm{d}y\textrm{d}x\right) \nonumber \\
 & = & \frac{\pi}{2\beta}-2\int_{0}^{R}\int_{0}^{\epsilon-\sqrt{R^{2}-x^{2}}}
	e^{-\beta\left(x^{2}+y^{2}\right)}\textrm{d}y\textrm{d}x
\end{eqnarray}
since the integral over $\mathcal{{D}}_{1}$ cancels. 

Thus consider two regimes for the distance $\epsilon$:
in the first, where $\epsilon \approx R$ (close to the boundary), we can make the approximation
$\exp\left({-\beta y^{2}}\right)\approx1$, since the distances $y$ from
the horizontal to the lower semi-circle in Fig. \ref{fig:annulus}
will be small, so we can approximate the integral in Eq. \ref{eq:neweq16}
\begin{eqnarray}\label{eq:diskmass}
\int_{0}^{R}\int_{0}^{\epsilon-\sqrt{R^{2}-x^{2}}}e^{-\beta\left(x^{2}+y^{2}\right)}\textrm{d}y\textrm{d}x & \approx & \int_{0}^{R}\int_{0}^{\epsilon-\sqrt{R^{2}-x^{2}}}e^{-\beta x^{2}}\textrm{d}y\textrm{d}x \nonumber \\
 & = & \frac{\sqrt{\pi}}{2\sqrt{\beta}}\epsilon-\int_{0}^{R}e^{-\beta x^{2}}\sqrt{R^{2}-x^{2}}\textrm{d}x\end{eqnarray}
such that
\begin{equation}
\mathcal{{M}}(\epsilon \approx R) = \frac{\pi}{2\beta}-\frac{1}{R\sqrt{\beta}}\left(\frac{\sqrt{\pi}}{4\beta}\right)+\left(R-\epsilon\right)\sqrt{\frac{\pi}{\beta}}+\mathcal{O}\left(\left(R-\epsilon\right)^2\right)
\end{equation}
after Taylor expanding Eq. \ref{eq:diskmass} for $\epsilon \approx R$, since it is in this regime that the main contribution to 
Eq. \ref{e:neweq6} comes from.

For the other regime (where $\epsilon \ll R$)
\begin{eqnarray}\label{eq:neweq18}
	\mathcal{{M}}(\epsilon \ll R) &\approx& \int_{0}^{2\pi}\int_{0}^{\infty}r' e^{-\beta r'^{2}} \textrm{d}r'\textrm{d}\theta \\
	&=& \pi / \beta
\end{eqnarray}
due to the exponential decay of the connectivity function, and so we have the probability of connection $P_{fc}$
\begin{eqnarray}\label{eq:4}
P_{fc} & \approx & 1 - \rho \int_{\mathcal{V}} e^{-\rho \int_{\mathcal{V}} \chi \left(x,y\right) H\left(r_{xy}\right) \mathrm{d}y } \mathrm{d}x \nonumber \\
 & = & 1 - \rho \int_{0}^{2\pi}\int_{0}^{L^{+}}\epsilon\exp\left(-\frac{\rho\pi}{\beta}\right)\textrm{d}\epsilon\textrm{d}\theta \nonumber \\
	   &  & - \rho\int_{L^{+}}^{R}\exp\left(-\rho\left(\left(\frac{\pi}{2\beta}-\frac{1}{R\sqrt{\beta}}\left(\frac{\sqrt{\pi}}{4\beta}\right)\right)+\left(R-\epsilon\right)\sqrt{\frac{\pi}{\beta}}\right)\right)\epsilon\textrm{d}\epsilon\textrm{d}\theta \nonumber \\
 & \approx & 1-\pi R^{2}\rho e^{-\frac{\rho\pi}{\beta}}-2\pi R\sqrt{\frac{\beta}{\pi}}e^{-\frac{\rho}{\beta}\left(\frac{\pi}{2}-\frac{1}{R\sqrt{\beta}}\left(\frac{\sqrt{\pi}}{4}\right)\right)}
\end{eqnarray}
where $L^+$ is the point where the two mass approximations equate. 
This approaches equation Eq. 38 of reference \cite{cef2012} as
$R\sqrt{\beta}\rightarrow\infty$, where the second term in the 
exponent of the final term in Eq. \ref{eq:4} is a `curvature correction' to the disk result in that report. 
Monte-Carlo simulations (where graphs are drawn from the graph-ensemble and enumerated should they connect), presented
in Fig. \ref{fig:montecarlo} alongside our approximation in Eq. \ref{eq:4}, corroborate our formula and show an improvement
on the result in \cite{cef2012}. The discrepancy at low density is expected since we only consider the
probability of a single isolated node.

We also highlight the interesting composition \ref{eq:4}. There is a bulk term 
(whose coefficient is proportional to the area of $\mathcal{D}$) and a boundary term (proportional to the circumference of $\mathcal{D}$).
This is discussed in greater detail in e.g. \cite{cef2012}, though we again emphasise the dominance
of the boundary term as $\rho \to \infty$.

\subsection{Small obstacles}
We now take the case where $r \ll r_0$ (but not necessarily zero), and take the outer perimeter $R \gg r_0$. 
We make the approximation that the small cone-like domain $\mathcal{{A}}_{c}$
(making up a portion
of the region visible $\mathcal{A}\left(x\right)$ to $x$ in the middle panel of
Fig. \ref{fig:trip}) is only significantly
contributing to the connectivity mass at small displacements $\epsilon$ from the obstacle, since at larger
displacements it thins and the
wedge-like region $\mathcal{A}\left(x\right) \setminus \mathcal{A}_{c}$ 
dominates.
Practically, it is $\mathcal{{A}}_{c}$ that presents the main integration difficulties,
so we approximate $H\left(r_{xy}\right)$ over this region
where the radial coordinate $r'\ll1$, using $\exp{-\beta r_{xy}'^{2}} \approx 1$ 
\begin{eqnarray}\label{eq:neweq17}
	\mathcal{M}\left( \epsilon \ll r_0 \right) & \approx & \int_{-\pi+\arcsin\left(\frac{r}{r+\epsilon}\right)}^{\pi-\arcsin\left(\frac{r}{r+\epsilon}\right)}\int_{0}^{\infty} e^{-\beta r'^{2}}r'\textrm{d}r'\textrm{d}\theta+2\int_{0}^{\arcsin\left(\frac{r}{r+\epsilon}\right)}\int_{0}^{\left(r+\epsilon\right)\cos(\theta)-\sqrt{r^{2}-\left(r+\epsilon\right)^{2}\sin^{2}(\theta)}}r'\textrm{d}r'\textrm{d}\theta \nonumber \\
 & = & \frac{1}{\beta}\left(\pi-\arcsin\left(\frac{r}{r+\epsilon}\right)\right)+\int_{0}^{\arcsin\left(\frac{r}{r+\epsilon}\right)}
	\left(\left(r+\epsilon\right)\cos\left(\theta\right)-\sqrt{r^{2}-\left(r+\epsilon\right)^{2}\sin^{2}\left(\theta\right)}\right)^2\textrm{d}\theta \nonumber \\
 & = & \frac{\pi}{\beta}+\left(r^{2}-\frac{1}{\beta}\right)\arcsin\left(\frac{r}{r+\epsilon}\right)+r\sqrt{2r\epsilon+\epsilon^{2}}-\frac{\pi}{2}r^{2}
\end{eqnarray}
For small $\epsilon$ (where the main contribution to Eq. \ref{e:neweq6} is found) we have
\begin{equation}
\mathcal{M}\left(\epsilon\ll r_0 \right) = \frac{\pi}{2\beta}+\frac{\sqrt{2}}{\beta\sqrt{r}}\epsilon^{1/2}+\frac{8\beta r^{2}-5}{6\beta\sqrt{2}r^{3/2}}\epsilon^{3/2}+\mathcal{O}\left(\epsilon^{2}\right)
\end{equation}
leaving us to integrate over the annulus
\begin{eqnarray} \label{eq:7}
P_{fc} & \approx & 1-\rho \int_{\mathcal{{A}}} e^{-\rho\mathcal{{M}}\left(x\right)}\textrm{d}x\nonumber \\
				   & \approx & 1 - \rho\int_{0}^{2\pi}\int_{0}^{L^{-}}\left(r+\epsilon\right)\exp\left(-\rho\left(\frac{\pi}{2\beta}+\frac{\sqrt{2}}{\beta\sqrt{r}}\epsilon^{1/2}+\frac{8\beta r^{2}-5}{6\beta\sqrt{2}r^{3/2}}\epsilon^{3/2}\right)\right)\textrm{d}\epsilon \textrm{d}\theta \nonumber \\ 
	& &  - \pi R^{2}\rho e^{-\frac{\rho\pi}{\beta}}-2\pi R\sqrt{\frac{\beta}{\pi}}e^{-\frac{\rho}{\beta}\left(\frac{\pi}{2}-\frac{1}{R\sqrt{\beta}}\left(\frac{\sqrt{\pi}}{4}\right)\right)} \nonumber \\
					& \approx & 1- 2\pi\rho\int_{0}^{L^{-}}\left(r+\epsilon\right) e^{-\frac{\rho\pi}{2\beta}}e^{-\rho\frac{\sqrt{2}}{\beta\sqrt{r}}\epsilon^{1/2}}\left(1-\rho\frac{8\beta r^{2}-5}{6\beta\sqrt{2}r^{3/2}}\epsilon^{3/2}\right)\textrm{d}\epsilon - \pi R^{2}\rho e^{-\frac{\rho\pi}{\beta}}-2\pi R\sqrt{\frac{\beta}{\pi}}e^{-\frac{\rho}{\beta}\left(\frac{\pi}{2}-\frac{1}{R\sqrt{\beta}}\left(\frac{\sqrt{\pi}}{4}\right)\right)} \nonumber \\
	                & \approx & 1 - \pi r^{2}\frac{2\beta^{2}}{\rho}e^{-\frac{\rho\pi}{2\beta}}-\pi R^{2}\rho e^{-\frac{\rho\pi}{\beta}}-2\pi R\sqrt{\frac{\beta}{\pi}}e^{-\frac{\rho}{\beta}\left(\frac{\pi}{2}-\frac{1}{R\sqrt{\beta}}\left(\frac{\sqrt{\pi}}{4}\right)\right)} 
 \end{eqnarray}
where $L^{-}$ is the point where the connectivity mass in the bulk meets our approximation $\mathcal{M}\left(\epsilon\ll r_0 \right)$ near the obstacle.
We numerically corroborate Eq. \ref{eq:7} in Fig. \ref{fig:montecarlo} using
Monte Carlo simulations. 

Note that this obstacle term is extremely small
compared to the other contributions in Eq. \ref{eq:7}, given its
coefficient decays linearly with $\rho$ and the factor of $\left(r\sqrt{\beta}\right)^2 \ll 1$. We conclude that 
a small internal perimeter of radius $r$ in any convex domain $\mathcal{V}$ results in
a negligible effect on connectivity. 

\subsection{Large obstacles}
For the large obstacle case ($r \gg r_0$)
\begin{eqnarray}\label{eq:8}
\mathcal{M}\left(\epsilon \ll r_0 \right) & \approx & 
	2\int_{0}^{\infty}\int_{0}^{\infty}e^{-\beta(x^{2}+y^{2})}\textrm{d}x\textrm{d}y+\int_{-\infty}^{\infty}\int_{0}^{
	\epsilon + \frac{1}{2r}x^2}e^{-\beta(x^{2}+y^{2})}\textrm{d}y\textrm{d}x\nonumber \\
 & = & \frac{\pi}{2\beta}+\frac{\sqrt{\pi}}
	{2\sqrt{\beta}}\int_{-\infty}^{\infty}e^{-\beta x^{2}}\textrm{erf}\left[\left(\epsilon+\frac{1}{2r}x^{2}\right)\sqrt{\beta}\right]\textrm{d}x\nonumber \\
\end{eqnarray}
yielding a power series in $\epsilon$
\begin{eqnarray}
\mathcal{M}\left(\epsilon \ll r_0 \right) & \approx & \frac{\pi}{2\beta}+\frac{\pi}{2\beta}\textrm{erf}\left[\sqrt{\beta}\epsilon\right]+\frac{1}{r\sqrt{\beta}}\left(\frac{\sqrt{\pi}}{4\beta}e^{-\beta\epsilon^{2}}\right)\nonumber \\
 & = & \frac{\pi}{2\beta}+\frac{1}{r\sqrt{\beta}}\left(\frac{\sqrt{\pi}}{4\beta}\right)+\frac{\sqrt{\pi}}{\sqrt{\beta}}\epsilon+\mathcal{O}\left(\epsilon^{3/2}\right)
\end{eqnarray}
This implies the connectivity mass is scaling in the same way as
for the outer boundary, but where the curvature correction (in the exponent of the last term in 
Eq. \ref{eq:7}) is of opposite sign. We therefore have
\begin{equation}\label{eq:11}
P_{fc} \approx 1-2\pi r\sqrt{\frac{\beta}{\pi}}e^{-\frac{\rho}{\beta}\left(\frac{\pi}{2}+\frac{1}{r\sqrt{\beta}}\left(\frac{\sqrt{\pi}}{4}\right)\right)}-
\pi R^{2}\rho e^{-\frac{\rho\pi}{\beta}}-2\pi R\sqrt{\frac{\beta}{\pi}}e^{-\frac{\rho}{\beta}\left(\frac{\pi}{2}-\frac{1}{R\sqrt{\beta}}\left(\frac{\sqrt{\pi}}{4}\right)\right)} 
\end{equation}
which is corroborated numerically in Fig. \ref{fig:montecarlo}.

This implies that large obstacles behave like separate, internal perimeters. In the large-domain limit (where the node numbers go
to infinity and the connection range is tiny compared to the large domain geometry), we can thus use
\begin{equation}\label{e:neweq15}
P_{fc} \approx 1-2\pi \left(R+r\right)\sqrt{\frac{\beta}{\pi}}e^{-\frac{\rho\pi}{2\beta}}-\pi(R^{2}-r^{2})\rho e^{-\frac{\rho\pi}{\beta}}
\end{equation}

\section{The spherical shell $\mathcal{{S}}$} \label{sec:four}

Consider now the spherical shell domain $\mathcal{{S}}$
of inner radius $r$ and outer radius $R$, which is the
three-dimensional analogue of the annulus. We again ask for the connection probability $P_{fc}$.

\subsection{Small spherical obstacles}
The region visible to the node at $x$ is again decomposed into two parts,
the three-dimensional version of $\mathcal{{A}}_{c}$,
called $\mathcal{{S}}_{c}$, and the rest of the region visible to $x$, denoted 
$\mathcal{{S}}(x)\setminus\mathcal{{S}}_{c}$. As in the annulus with the small obstacle, we
approximate $H\left(r_{xy}\right)$ over this region
where the radial coordinate $r'\ll1$ (which holds for $\epsilon \ll 1$ where the main contribution to the connectivity mass
is found), using $\exp{-\beta r_{xy}'^{2}} \approx 1$ such
that the contribution to the connectivity mass over $\mathcal{M}_{\mathcal{S}_{c}}\left(\epsilon\right)$
is
\begin{eqnarray}
\mathcal{M}_{\mathcal{S}_{c}}\left(\epsilon\right) & = & \int_{\mathcal{{S}}_{c}}r'^{2}e^{-\beta r'^{2}}\sin\theta \textrm{d}r'\textrm{d}\theta \textrm{d}\varphi \nonumber \\
 & \approx & \int_{\mathcal{{S}}_{c}}r'^{2}\sin\theta \textrm{d}r'\textrm{d}\theta \textrm{d}\varphi
\end{eqnarray}
We evaluate this by breaking up $\mathcal{{S}}_{c}$ into the area of a cone of
radius $\lambda$, height $h$ and apex angle $2\theta_{c}$
\begin{eqnarray}
\lambda&=&\frac{r}{r+\epsilon}\sqrt{2r\epsilon+\epsilon^{2}}\\
h&=&\frac{2r\epsilon+\epsilon^{2}}{r+\epsilon}\\
\theta_{c}&=&\arcsin\left(\frac{r}{r+\epsilon}\right)
\end{eqnarray}
(with the apex at a distance $\epsilon$ from the obstacle),
and the spherical segment (which on removal from the cone creates
the shape of $\mathcal{S}_{c}$)
\begin{eqnarray}
\mathcal{{M}}_{\mathcal{{S_{C}}}}\left(\epsilon\right) & = & \frac{1}{3}\pi \lambda^{2}h-\frac{1}{6}\pi\left(r+\epsilon-h\right)\left(3\lambda^{2}+\left(r+\epsilon-h\right)^{2}\right)\nonumber \\
 & = & \frac{\epsilon^{2}\pi r^{2}\left(\epsilon+2r\right)^{2}}{3\left(\epsilon+r\right)^{3}}-\frac{\epsilon^{2}\pi r^{3}\left(2\epsilon+3r\right)}{3\left(\epsilon+r\right)^{3}}\nonumber \\
 & = & \frac{\epsilon^{2}\pi r^{2}}{3\left(\epsilon+r\right)}
\end{eqnarray}
Adding the mass over $\mathcal{{S}}\left( x \right)\setminus\mathcal{{S}}_{c}$, we use the fact that
the full solid angle available to a bulk node is $4\pi$ and that the angle $\omega\leq\Omega$
available to the node at $x$ is
\begin{eqnarray}
	\omega & = & \frac{1}{4\pi}\int_{0}^{2\pi}\int_{0}^{\theta_{c}}\sin\left(\theta\right)\textrm{d}\theta \textrm{d}\varphi \nonumber \\
 & = & \frac{1}{2}\left(1-\cos\left(\arcsin\left(\frac{r}{r+\epsilon}\right)\right)\right)\nonumber \\
 & = & \frac{1}{2}\left(1-\sqrt{1-\left(\frac{r}{r+\epsilon}\right)^{2}}\right)
\end{eqnarray}
such that
\begin{equation}
\int_{\mathcal{{S}}\left( x \right)\setminus\mathcal{{S}}_{c}}r'^{2}e^{-\beta r'^{2}}\sin\theta \textrm{d}r'\textrm{d}\theta \textrm{d}\varphi =\frac{\pi\sqrt{\pi}}{\beta\sqrt{\beta}}\left(1-\frac{1-\sqrt{1-\left(\frac{r}{r+\epsilon}\right)^{2}}}{2}\right)
\end{equation}
We then have $\mathcal{{M}}\left( \epsilon \ll r_0 \right)$
\begin{eqnarray}
\mathcal{{M}}\left( \epsilon \ll r_0 \right)&\approx&\frac{\epsilon^{2}\pi r^{2}}{3\left(\epsilon+r\right)}+\frac{\pi\sqrt{\pi}}{\beta\sqrt{\beta}}\left(1-\frac{1-\sqrt{1-\left(\frac{r}{r+\epsilon}\right)^{2}}}{2}\right)\\
&=&\frac{\pi\sqrt{\pi}}{2\beta\sqrt{\beta}}+\frac{\pi\sqrt{\pi}}{\beta\sqrt{\beta}}\frac{1}{\sqrt{2r}}\epsilon^{1/2}+\frac{3\pi^{3/2}}{4\sqrt{2}\left(r\beta\right)^{3/2}}\epsilon^{3/2}+\mathcal{O}\left(\epsilon^{2}\right)\end{eqnarray}
which implies that small spherical obstacles reduce the connection probability
within the unobstructed sphere domain $\mathcal{S}_{r=0}$ to give a connection probability of
\begin{eqnarray}
	P_{fc}^{\mathcal{S}_{r \ll r_0}} & \approx & P_{fc}^{\mathcal{S}_{r=0}} - \rho e^{-\rho\left(\frac{\pi\sqrt{\pi}}{2\beta\sqrt{\beta}}\right)}\int_{0}^{2\pi}\int_{0}^{\pi}\int_{0}^{L_{\mathcal{{S}}}^{-}}\left(r+\epsilon\right)^{2}\sin\left(\theta\right)e^{-\rho\left(\frac{\pi\sqrt{\pi}}{2\beta\sqrt{\beta}}+\frac{\pi\sqrt{\pi}}{\beta\sqrt{\beta}}\frac{1}{\sqrt{2r}}\epsilon^{1/2}+\frac{3\pi^{3/2}}{4\sqrt{2}\left(r\beta\right)^{3/2}}\epsilon^{3/2}\right)}\textrm{d}\epsilon \textrm{d}\theta \textrm{d}\varphi \nonumber \\
																			& \approx & P_{fc}^{\mathcal{S}_{r=0}} - \frac{4}{3}\pi r^{3}\left(\frac{12\beta^{3}}{\rho\pi^{3}}\right)e^{-\rho\left(\frac{\pi\sqrt{\pi}}{2\beta\sqrt{\beta}}\right)}
\end{eqnarray}\newline

\begin{figure}
\noindent
\begin{centering}
	\includegraphics[scale=0.22]{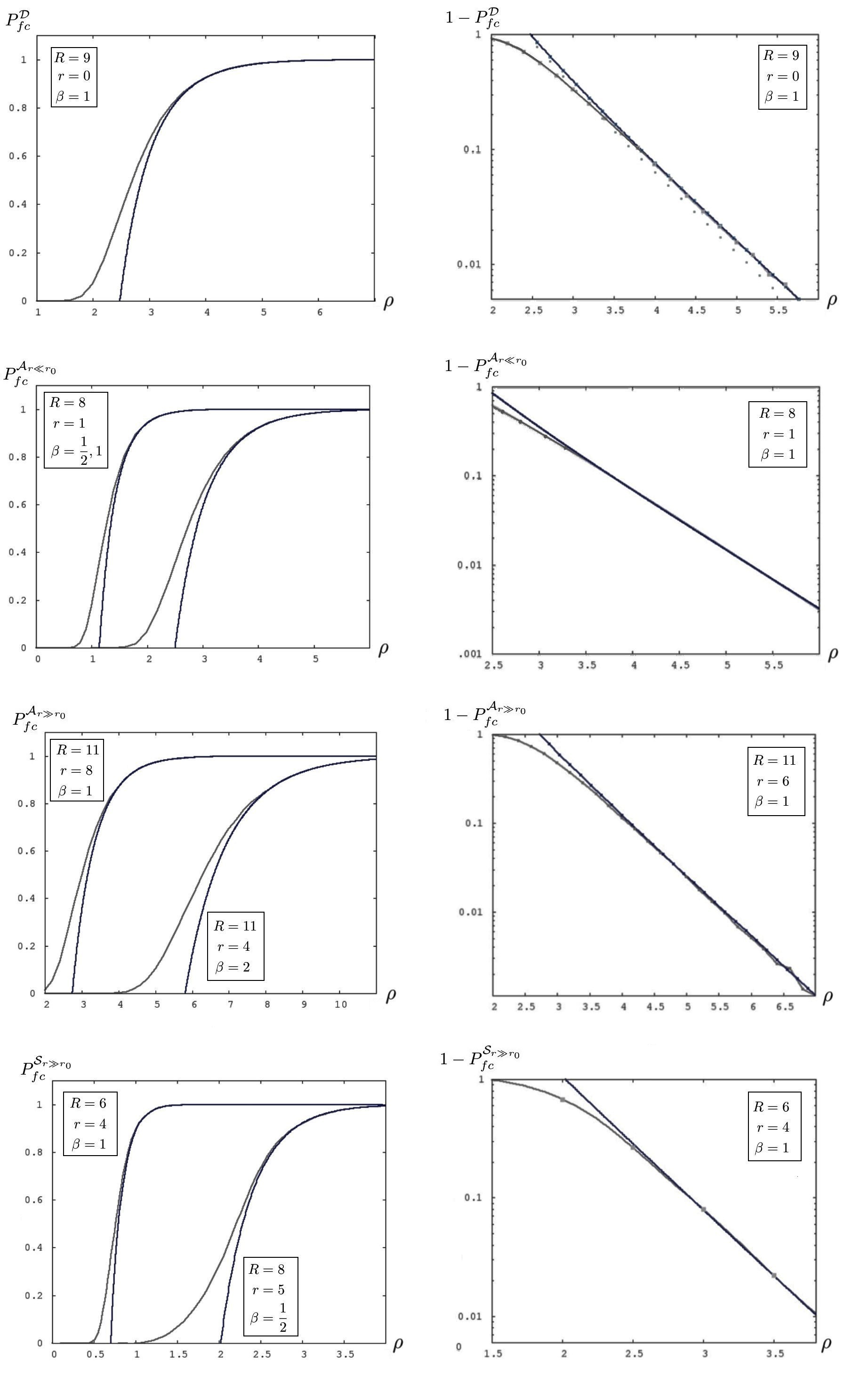}
\par\end{centering}
\caption{(Colour online) Monte Carlo simulations: We use Monte Carlo methods to estimate the connection probability of 
	soft random geometric graphs drawn inside various annuli and spherical shells $\mathcal{A}$ and $\mathcal{S}$ respectively. Every curve is
	compared with our analytic predictions (darker line) from Eqs. \ref{eq:4}, \ref{eq:7}, \ref{eq:11} and \ref{eq:10} (where indicated). 
	The discrepancy at low density is expected due to the fact we calculate only the probability of a single isolated node, given the results
	in e.g. \cite{penrose2013}.}\label{fig:montecarlo}
\end{figure}

\subsection{Large spherical obstacles}
For large obstacles ($r \gg r_0$), we extend Eq. \ref{eq:8} into the third dimension. 
$\mathcal{M}\left( \epsilon \ll r_0 \right)$ thus becomes
\begin{equation}
\mathcal{M}\left( \epsilon \ll r_0 \right) \approx 
4\int_{0}^{\infty}\int_{0}^{\infty}\int_{0}^{\infty}\textrm{d}x\textrm{d}y\textrm{d}ze^{-\beta(x^{2}+y^{2}+z^{2})}+\int_{-\infty}^{\infty}\int_{-\infty}^{\infty}\int_{0}^{\nu\left(x,z\right)}\textrm{d}y\textrm{d}x\textrm{d}ze^{-\beta(x^{2}+y^{2}+z^{2})}
\end{equation}
where $\nu\left(x,z\right)=\epsilon+\frac{1}{2r}\left(x^{2}+z^{2}\right)$, yielding
\begin{eqnarray}\label{eq:spheremass}
\mathcal{M}\left( \epsilon \ll r_0 \right) & \approx & \frac{\pi\sqrt{\pi}}{2\beta\sqrt{\beta}}+\frac{\pi\left(r\sqrt{\beta\pi}\textrm{erf}\left[\epsilon\sqrt{\beta}\right]+e^{-\beta\epsilon^{2}}\right)}{2r\beta^{2}}\nonumber \\
 & = & \frac{\pi\sqrt{\pi}}{2\beta\sqrt{\beta}}+\frac{\pi}{2\beta^{2}r}+\frac{\pi}{\beta}\epsilon+\mathcal{O}\left(\epsilon^{2}\right)
\end{eqnarray}
implying the connection probability is
\begin{eqnarray}
P_{fc}^{\mathcal{S}_{r \gg r_0}} & \approx & P_{fc}^{\mathcal{S}_{r=0}} - \rho e^{-\rho\left(\frac{\pi\sqrt{\pi}}{2\beta\sqrt{\beta}}+\frac{\pi}{2\beta^{2}r}\right)}\int_{0}^{2\pi}\int_{0}^{\pi}\int_{0}^{L_{\mathcal{{S}}}^{-}}\left(r+\epsilon\right)^{2}\sin\left(\theta\right)\textrm{d}\epsilon \textrm{d}\theta \textrm{d}\varphi e^{-\rho\left(\frac{\pi}{\beta}\epsilon\right)}\nonumber \\
								 &\approx&  P_{fc}^{\mathcal{S}_{r=0}} - 4\pi r^{2}\left(\frac{\beta}{\pi}\right)e^{-\rho\left(\frac{\pi\sqrt{\pi}}{2\beta\sqrt{\beta}}+\frac{1}{R\sqrt{\beta}}\left(\frac{\pi}{2\beta\sqrt{\beta}}\right)\right)}
\end{eqnarray}\newline
where $L^-$ is the point where our mass approximation in Eq. \ref{eq:spheremass} is equal to the mass
in the bulk of the sphere $\left(\pi/\beta \right)^{3/2}$
(given the argument used for the two-dimensional case in Eq. \ref{eq:neweq18}).

We now have the connection probability in the spherical shell $\mathcal{S}$
\begin{multline}\label{eq:10}
	P_{fc}^{\mathcal{S}} \approx 
1-\frac{4\pi}{3}\left(R^{3}-r^{3}\right)\rho e^{-\rho\left(\frac{\pi\sqrt{\pi}}{\beta\sqrt{\beta}}\right)}
-4\pi R^{2}\left(\frac{\beta}{\pi}\right)e^{-\rho\left(\frac{\pi\sqrt{\pi}}{2\beta\sqrt{\beta}}-\frac{1}{R\sqrt{\beta}}
\left(\frac{\pi}{2\beta\sqrt{\beta}}\right)\right)}\\
-	\begin{cases}
		\frac{4}{3}\pi r^{3}\left(\frac{12\beta^{3}}{\rho\pi^{3}}\right)e^{-\rho\left(\frac{\pi\sqrt{\pi}}{2\beta\sqrt{\beta}}\right)}  & \quad \text{if } r \ll r_0 \\
		4\pi r^{2}\left(\frac{\beta}{\pi}\right)e^{-\rho\left(\frac{\pi\sqrt{\pi}}{2\beta\sqrt{\beta}}
-\frac{1}{R\sqrt{\beta}}\left(\frac{\pi}{2\beta\sqrt{\beta}}\right)\right)}  & \quad \text{if } r \gg r_0
	\end{cases}
\end{multline}
which is corroborated numerically in Fig. \ref{fig:montecarlo} (but only for the large obstacle case, given the negligable
size of the small obstacle term in comparison to the bulk and boundary contributions). Just as with the annulus, 
small spherical obstacles thus
have little impact on connectivity, and large spherical obstacles behave like separate perimeters. This
behaviour is likely the same for all dimensions $d>3$, where the geometry is a hypersphere containing a convex
$d$-dimensional obstacle (which one might call a hyper-annulus).

\section{Scenarios where obstacle effects are dominant} \label{sec:five}

In the previous sections, we have provided approximations for the probability of a single isolated node appearing
within both the annulus $\mathcal{A}$ and spherical shell $\mathcal{S}$. We draw three main conclusions:

\begin{enumerate}
	\item Small obstacles holes in the domain have little effect on connectivity (in all dimensions $d \geq 2$).\\
	\item Large obstacles holes disrupt connectivity as separate domain boundaries (in all dimensions $d \geq 2$). \\
	\item As $\rho \to \infty$, the effect of a convex obstacle of any size is quickly dominated by that of the
		enclosing perimeter (in all dimensions $d \geq 2$).
\end{enumerate}

One may therefore be forgiven for suggesting that obstacles have little impact on connectivity in dense networks.
This is not always true, and so in the next few subsections we highlight important situations where these convex obstructions are
\textit{essential} to connectivity, since the `small correction term' provided by obstacle analysis in the dense
network limit becomes (in some parameter regime) the dominant contribution to $P_{fc}$.

\subsection{Multiple convex obstacles distributed over $\mathcal{V}$}
\begin{figure}
\noindent
\begin{centering}
	\includegraphics[scale=0.5]{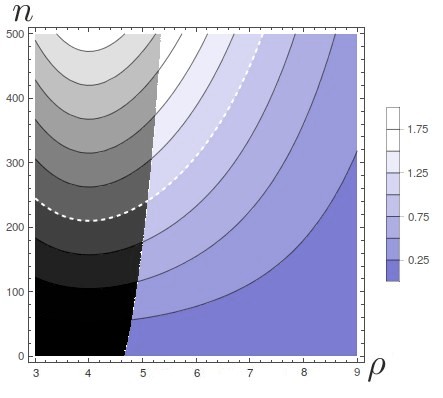} \hspace{10mm} \includegraphics[scale=0.5]{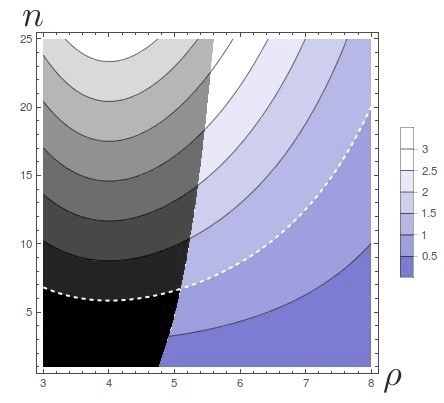}
\par\end{centering}
\caption{(Colour online) Many obstacles: Taking the Sinai-like domain in Fig. \ref{fig:annulus}
	with side $L=100$ and containing $n$ circular obstacles of small radii $r=1$
	(\textit{left-hand phase diagram})	and large radii $r=6$ (\textit{right-hand phase diagram}),
we plot the ratio of the first term in Eq. \ref{eq:manyobs2} with the sum of the components
in Eq. \ref{eq:manyobs1}, whose magnitude is indicated by the color gradient on the respective legend. The obstacle effects dominate when
this exceeds unity, highlighted by the dotted line on each graph. Also, the regions where
$P_{fc}$ is (predicited by our formulas to be) below $4/5$ are faded to gray tones, indicating regions where our approximations to the
connection probability begin to lose their accuracy. $\beta = 1$ throughout.}
\label{fig:phase}
\end{figure}
Given that the holes are not too close, their effects add up in a linear fashion such that they potentially outweigh
the effect of the boundary. To highlight this, take the Sinai-like domain 
in the right hand panel of Fig. \ref{fig:annulus}. Without obstacles, we have

\begin{eqnarray}\label{eq:manyobs1}
	P_{fc} = 1 - L^{2}\rho e^{-\frac{\pi}{\beta} \rho } - 4 L\sqrt{\frac{\beta}{\pi}}  e^{-\frac{\pi}{2\beta} \rho} - \frac{16\beta}{\rho\pi}  e^{-\frac{\pi}{4\beta} \rho}
\end{eqnarray}
taken from \cite{cef2012}. This is composed of a bulk term, a boundary term and a corner term. As we have seen, 
introducing $n$ circular obstacles of various radii $r_i$ will reduce this connection probability such that we have
\begin{eqnarray}\label{eq:manyobs2}
	P_{fc} = 1 - \sum_{i=1}^{n}\pi r_{i}^{2} \left(\frac{2\beta^{2}}{\rho}\right)e^{-\frac{\rho\pi}{2\beta}}
	- \left(L^{2}-\sum_{i=1}^{n}\pi r_{i}^{2}\right)\rho e^{-\frac{\pi}{\beta} \rho } - 4L \sqrt{\frac{\beta}{\pi}}  e^{-\frac{\pi}{2\beta} \rho} - \frac{16\beta}{\rho\pi}  e^{-\frac{\pi}{4\beta} \rho}
\end{eqnarray}
which holds whenever the obstacles are separated from each other and the boundary 
by at least $2r_0$. Fig.
\ref{fig:phase} presents two phase plot that demonstrate how the obstacle effects can become
dominant given a certain number of obstacles $n$. As we pass through
the moderate density regime, the obstacle
effects pass through a phase of significance greater than the sum of the rest of the geometric 
contributions to $P_{fc}$ (i.e. the bulk, square
perimeter and four corners).

\subsection{Surfaces without boundary}

Boundary effects can be removed by working on surfaces without an enclosing perimeter.
Examples include the flat torus (popular in
rigorous studies but difficult to realise in wireless networks), and the sphere.
Thus as $\rho \to \infty$ the obstacle effects are
the dominant contribution to $P_{fc}$.

This may be of interest to pure mathematicians studying random graphs for purposes
outside communication theory \cite{penrosebook}. Fractal obstacles may be of particular interest \cite{dettmann2014}.

\subsection{Quasi-1D regime $r \approx R$}

Note that as the width of the annulus goes to zero, the approximation used in
Eq. \ref{e:neweq6} (that connectivity is the same as no isolated nodes) breaks down.
The graph now disconnects by forming two clusters separated from each other by two
unpopulated strips of width usually greater than $ r_0$. We
call this situation (where $R \approx r$) the `quasi-1D' regime, deferring its treatment to a later study. We 
emphasise that one-dimensional random geometric graphs are particularly interesting, since they provide a test-bed
for other theories that may be difficult to study initially in dimensions $d \geq 2$.

\section{Conclusions} \label{sec:six}

We have derived semi-rigorous analytic formulas for the connection probability 
of soft random geometric graphs drawn inside various annuli and shells 
(of inner radius $r$ and outer radius $R$) given the link formation probability between two nodes
is an exponentially decaying function of their Euclidean separation. This models the 
Rayleigh fading of radio signal propagation within a wireless \textit{ad hoc} network.

We have thus extended the soft connection model into 
simple non-convex spaces based on
circular or spherical obstacles (rather than fractal boundaries \cite{dettmann2014}, 
internal walls \cite{orestis2013} or fixed obstacles on a grid \cite{almiron2013}). We
highlight situations where obstacles are (and are not) important influences on connectivity:

\begin{enumerate}

	\item Small obstacles have little impact on connectivity. \\
	\item Large obstacles have a similar impact on connectivity as the enclosing perimeter, but their effects are
		dominated by the boundary as $\rho \to \infty$.
	\item Multiple obstacles can have the dominant effect on connection within density regimes that are significant for
		various areas
		of application, particularly \textit{ad hoc} communication networks deployed in urban environments. 
		5G wireless networks are an example of this scenario.

\end{enumerate}
Further topics of study include the quasi one-dimensional regime, where connectivity is not governed by isolated nodes.\\

Understanding the connectivity of these spatially embedded graphs in non-convex domains is a crucial
enabler for the reality of 5G wireless networks, particularly if these multi-hop relay systems form
in cluttered, urban environments (which is likely). Limiting scenarios (such as `many obstacles' and the `quasi-1D' regime) prove
to be particularly interesting.

\subsection*{\em Acknowledgements} The authors wish to thank Suzanne Binding and the directors of the Centre for Doctoral Training in Communications
at the University of Bristol, alongside the directors of the Toshiba Telecommunications Research Laboratory (Europe) for their continued support. 
Thanks also go to the anonymous referee who suggested important clarifications to section \ref{sec:two}, and to the other referee for their very useful comments.
They also thank Justin Coon, Tom Kealy, Leo Laughlin, Jon Keating and David Simmons for many helpful discussions.


\begin{thebibliography}{10}

\bibitem{penrose2013}
M.~D. Penrose, ``Connectivity of soft random geometric graphs,'' {\em to appear
  in Ann. Appl. Probab., available at arXiv:1311.3897}, 2013.

\bibitem{newmanbook}
M.~E.~J. Newman, {\em Networks: An Introduction}.
\newblock New York, NY, USA: Oxford University Press, 2010.

\bibitem{brecht2013}
J.~von Brecht, T.~Kolokolnikov, A.~Bertozzi, and H.~Sun, ``Swarming on random
  graphs,'' {\em J. Stat. Phys.}, vol.~151, no.~1-2, pp.~150--173, 2013.

\bibitem{eubank2004}
S.~Eubank, H.~Guclu, V.~S.~A. Kumar, M.~V. Marathe, A.~Srinivasan,
  Z.~Toroczkai, and N.~Wang, ``Modelling disease outbreaks in realistic urban
  social networks,'' {\em Nature}, vol.~429, no.~6988, pp.~180--184, 2004.

\bibitem{amin2013}
M.~Amin, ``Energy: The smart grid solution,'' {\em Nature}, vol.~499, no.~7457,
  pp.~145--147, 2013.

\bibitem{tsebook}
D.~Tse and P.~Viswanath, {\em Fundamentals of Wireless Communication}.
\newblock Cambridge University Press, 2005.

\bibitem{mao2011}
G.~Mao and B.~D. Anderson, ``{On the Asymptotic Connectivity of Random Networks
  under the Random Connection Model},'' {\em {INFOCOM, Shanghai, China}},
  p.~631, 2011.

\bibitem{orestis2014}
O.~Georgiou, C.~Dettmann, and J.~Coon, ``Network connectivity: Stochastic vs.
  deterministic wireless channels,'' {\em Proc. IEEE ICC 2014, Sydney,
  Australia}, pp.~77--82, 2014.

\bibitem{hagenggi2009}
M.~Haenggi, J.~Andrews, F.~Baccelli, O.~Dousse, and M.~Franceschetti,
  ``Stochastic geometry and random graphs for the analysis and design of
  wireless networks,'' {\em Selected Areas in Communications, IEEE Journal on},
  vol.~7, no.~27, pp.~1029--1046, 2009.

\bibitem{li2009}
J.~Li, L.~Andrew, C.~Foh, M.~Zukerman, and H.~Chen, ``Connectivity, coverage
  and placement in wireless sensor networks,'' {\em Sensors}, vol.~9,
  pp.~7664--7693, 2009.

\bibitem{tohbook}
C.~K. Toh, {\em Ad Hoc Mobile Wireless Networks: Protocols and Systems}.
\newblock Prentice Hall, 2001.

\bibitem{halldorsson2014}
M.~M. Halld\'{o}rsson and T.~Tonoyan, ``How well can graphs represent wireless
  interference?,'' {\em Proc. Forty-Seventh Annual ACM Symposium on the Theory of Computing STOC '15}, pp.~635--644, 2015.

\bibitem{cef2012}
J.~Coon, C.~Dettmann, and O.~Georgiou, ``Full connectivity: Corners, edges and
  faces,'' {\em J. Stat. Phys.}, vol.~147, no.~4, pp.~758--778, 2012.

\bibitem{orestis2015}
O.~Georgiou, C.~Dettmann, and J.~Coon, ``Connectivity of networks with general
  connection functions,'' {\em preprint available at arXiv:1411.3617}, 2014.

\bibitem{clark1991}
B.~Clark, C.~Colbourn, and D.~Johnson, ``Unit disk graphs,'' {\em Discrete
  Mathematics}, vol.~86, no.~1--3, pp.~165--177, 1991.

\bibitem{erdos19602}
P.~Erd\H{o}s and A.~R\'{e}nyi, ``On random graphs,'' in {\em {Publ. Math.
  Debrecen}}, vol.~6, pp.~290--297, 1959.

\bibitem{gilbert1961}
G.~Gilbert, ``Random plane networks,'' {\em SIAM J.}, vol.~9, no.~4,
  pp.~533--543, 1961.

\bibitem{orestis2013}
O.~Georgiou, C.~Dettmann, and J.~Coon, ``Network connectivity through small
  openings,'' {\em Proc. ISWCS '13, Ilmenau, Germany}, pp.~1--5, 2013.

\bibitem{kog105}
A.~P. Giles, O.~Georgiou, and C.~P. Dettmann, ``Betweenness centrality in dense
  random geometric networks,'' {\em Proc. IEEE ICC 2015, London, UK}, 2015.

\bibitem{almiron2013}
M.~G. Almiron, O.~Goussevskaia, A.~A. Loureiro, and J.~Rolim, ``Connectivity in
  obstructed wireless networks: From geometry to percolation,'' in {\em
  Proceedings of the Fourteenth ACM International Symposium on Mobile Ad Hoc
  Networking and Computing, Bangalore, India}, pp.~157--166, 2013.

\bibitem{bocus2013}
 O.~Georgiou, M.~Z. Bocus, M.~R. Rahman, C.~P. Dettmann and J.~P. Coon, ``Keyhole and
  reflection effects in network connectivity analysis,'' {\em IEEE Commun. Lett.} vol. 19, no. 3, pp.~427--430, 2015.

\bibitem{dettmann2014}
C.~P. Dettmann, O.~Georgiou, and J.~P. Coon, ``More is less: Connectivity in
  fractal regions,'' {\em Proc. IEEE ICC 2015, London, UK}, 2015.

\bibitem{coon2014}
J.~P. Coon, O.~Georgiou, and C.~P. Dettmann, ``Connectivity in dense networks
  confined within right prisms,'' {\em 12th International Symposium on Modeling and Optimization in Mobile, Ad Hoc, and Wireless Networks, WiOpt 2014,}
  Hammamet, Tunisia, 2014. 

\bibitem{kingmanbook}
J.~F.~C. Kingman, {\em Poisson Processes}.
\newblock Oxford University Press, 1993.

\bibitem{waltersreview}
M.~Walters, ``{Random Geometric Graphs},'' in {\em {Surveys in Combinatorics
  2011}} ({Robin Chapman}, ed.), Cambridge University Press, 2011.

\bibitem{gupta1998}
P.~Gupta and P.~R. Kumar, ``Critical power for asymptotic connectivity,'' in
  {\em Proc. 37th IEEE Conference on Decision and Control}, Tampa, Florida,
  vol. 1, pp.~1106--1110, 1998.

\bibitem{penrosebook}
M.~D. Penrose, {\em Random Geometric Graphs}.
\newblock Oxford University Press, 2003.

\bibitem{coon2012}
J.~Coon, C.~Dettmann, and O.~Georgiou, ``Impact of boundaries on fully
connected random geometric networks,'' {\em Phys. Rev. E}, vol.~85, 011138, 2012.

\bibitem{sklar1997}
B.~Sklar, ``Rayleigh fading channels in mobile digital communication systems.
  Part I: Characterization,'' {\em IEEE Communications Magazine}, vol.~35,
  no.~7, 1997.

\end{thebibliography}
\end{document}